\documentclass[11pt]{article}

\usepackage{amsmath,amsthm,verbatim,amssymb,amsfonts,amscd, graphicx}
\usepackage{graphics}
\usepackage{hyperref}
\usepackage{authblk}
\usepackage{lscape}
\usepackage{multirow}
\usepackage{float}
\theoremstyle{plain}
\numberwithin{equation}{section}
\newtheorem{theorem}{Theorem}[section]

\newtheorem{corollary}{Corollary}[theorem]

\newtheorem{definition}{Definition}[section]
\date{}
\title{\textbf{Controllability Analysis of Motion of Artificial Satellite Under the Effect of Oblateness of the Earth}}
\author[1]{Jaita Sharma\thanks{jaita.sharma-appmath@msubaroda.ac.in}}
\author[1]{B. S. Ratanpal\thanks{bharatratanpal@gmail.com}}
\author[1]{Shivam Munshi\thanks{munshishivamdev120@gmail.com}}
\author[1]{Vishant Shah\thanks{vishantmsu83@gmail.com}}

\affil[1]{\small{Department of Applied Mathematics, Faculty of Technology \& Engineering, The Maharaja Sayajirao University of Baroda, Vadodara - 390 001, India}}

\begin{document}
\maketitle
\begin{abstract}
In this article we have studied the controllability of artificial satellite under the effect of zonal harmonic $J_2$ in cylindrical polar coordinates systems. Seven different cases of thrusters in various directions have been analyzed and it is found that the system is controllable if we apply thrusters in either $r$, $\theta$ and $z$ or $\theta$ and $z$ direction. The equations governing motion of satellite have been linearized and Kalman controllability test is applied to check the controllability of the system. We have also derived controller $u$ for the linearized system. The trajectory of the system have been plotted to show the controllability of the system.
\end{abstract}


\section{Introduction}
\label{sec:1}
Artificial satellites play very important role in navigation, communication, monitoring environment of the earth etc. \cite{Awange,Palancz}.
Many researchers studied motion of artificial satellite using analytic, semi-analytic and numerical methods. King-Hele\cite{King}
solved two-body problem of satellite, analytically by considering oblateness of Earth. Raj \cite{Raj} regularized equation of motion by applying KS transformations \cite{SS} and solved these equations of motion by considering atmospheric drag. 
Sehnal \cite{Sehnal} studied the motion of artificial satellite by considering perturbation due to upper terrestrial atmosphere. Knowles \textit{et.al.} \cite{KPTN} analyze the sample orbit from sensor data as well as orbital elements, during the period 14 July 2000, they found that geomagnetic storms driven by solar eruption have significant effect on the total density of the upper atmosphere in the altitude range 250 -- 1000 k.m., which causes a measurable effect on the orbit of resident space object.  
Yan and Kapila \cite{YK} developed the dynamical equations of satellite motion around oblate earth using spherical rotating frame and using this dynamics they derived conditions under which osculating plane of motion of satellite remains fixed. Khalil \cite{Khalil} developed analytical solution by considering atmospheric drag and oblateness of earth up to 4th order zonal harmonic using Hamiltonian mechanics. Bezdv\v{e}k and Vokrouhlick\'{y} \cite{BV} presented a semi-analytic theory for small eccentric orbit by considering oblateness of earth up to 9th order zonal harmonic of the earth and atmospheric drag. In this they considered empirical model TD88 of the neutral atmosphere density distribution for atmospheric drag. They also compared their predictions with the orbital data of several real-world artificial satellites. Hassan \textit{et. al.} \cite{HHB} regularized equations of perturbed motion due to oblateness of Earth using KS transformations and derived algorithm to solve these equations using Picard's method. Chen and Jing \cite{CJ} studied relative motion of satellite under the effect of the oblateness of earth and atmospheric drag.
Using Lie group variational approach Lee \textit{et. al.} \cite{LSSC} simulated rotational dynamics of satellite. Formation flight of artificial satellite under the effect of aerodynamics forces was studied by Reid and Misra \cite{RM}.
Xu and Chen \cite{XTCY} derived analytical solution in terms of Keplerian angular elements of satellite orbit under effect of atmospheric drag. Effect on the orbit of satellite Cosmos1484 under the effect of earth oblateness and atmospheric drag have been studied by Al-Bermani \textit{et. al.} \cite{AAAB}.
Using Lie transformations, Delhaise \cite{Delhaise} derived analytical solution of motion of satellite by considering gravity and air drag.

Sharma \textit{et.al.} \cite{SRPSCB,SRPS} studied the motion of satellite with different initial velocities and computed orbital elements by considering oblateness of earth and combined effect of the oblateness of Earth and atmospheric drag. They have also computed the time at which satellite will hit the Earth. To have satellite in correct orbit for longer time it is necessary to put controller that controls the motion.  
 Hajovsky \cite{Hajo} used atmospheric drag as a controller to control the trajectory of artificial satellite. B. Palancz \cite{Palancz, Palancz1} used pole placement to control trajectory of the artificial satellite. Recently Lamba \cite{Lamba}, discussed controllability, observability and stability problem concerned with artificial satellite using state space method. However he took two dimensional model which leads to sets of four equations in polar form. 

In this work we consider motion of satellite under the effect of $J_2$ zonal harmonic in cylindrical polar coordinate system and studied controllability of motion by plugging controllers (in form of thrusters) in various directions. It has been observed that the motion of satellite is controllable if controllers are kept in $r$, $\theta$ \& $z$ directions and $\theta$ \& $z$ directions. We also studied trajectory controllability of satellite.

\section{Preliminaries}
\label{sec:2}
In real life, most of the systems are nonlinear in nature and this nonlinearity creates difficulty in finding solution of the system. Hence it is required to approximate the nonlinear system by the appropriate linear system.

The motion of artificial satellite under the effect of zonal harmonic $J_2$ is modelled in terms of system of nonlinear differential equations. Here, we introduce the concept of linear control theory followed by linearization of nonlinear control systems \cite{Brockett}. 
\subsection{Linear Control Systems}
Consider linear control system,
\begin{eqnarray}
 \nonumber \dot{x}(t)&=&A(t)x(t)+B(t)u(t),\\
 x(t_0)&=&x_0, \label{lcs1}
\end{eqnarray}
where, $x_0,x(t)\in \mathbb{R}^n$ for all $t\in [t_0,t_1]$, $u\in L^2([t_0,t_1],\mathbb{R}^m)$. The matrices $A(t)$ and $B(t)$ are of order $n\times n$ and $n\times m$ respectively.\\
Let $\Phi(t,t_0)$ be the transition matrix of the homogeneous system $\dot{x}(t)=A(t)x(t)$ with initial condition $x(t_0)=x_0$ then solution of the system \eqref{lcs1} is given by,
\begin{equation}
 x(t)=\Phi(t,t_0)x_0+\int^{t}_{t_0}\Phi(t,s)B(s)u(s)ds\label{lcs2}.
\end{equation}
\begin{definition}
 The system \eqref{lcs1} is controllable over the interval $[t_0,t_1]$, if each pair of vectors $x_0$ and $x_1$ in $\mathbb{R}^n$ there is a control $u\in L^2([t_0,t_1],\mathbb{R}^m)$ such that the solution of \eqref{lcs1} satisfies $x(t_1)=x_1$. This means there is a control $u$ satisfying 
 \begin{equation*}
  x_1=\Phi(t_1,t_0)x_0+\int^{t_1}_{t_0} \Phi(t_1,s)B(s)u(s)ds.
 \end{equation*}
\end{definition}
 \begin{theorem}
  The system \eqref{lcs1} is controllable if and only if the controllability grammian of the system defined by $W(t_0,t_1)=\int^{t_1}_{t_0}\Phi(t_1,s)B(s)B^*(s)\Phi^*(t_1,s)ds$ is invertible and control $u$ of the system \eqref{lcs1} is given by 
  \begin{equation*}
   u(t)=B^*(t)\Phi^*(t_1,t)W^{-1}(t_0,t_1)[x_1-\Phi(t_1,t_0)]. \label{lcs3}
  \end{equation*}
  \end{theorem}
However if the system is time invariant, conditions reduces to Kalmann condition which is given by,
\begin{corollary}
If matrices $A$ and $B$ are two time invariant matrices of the system \eqref{lcs1} then the system is controllable if and only if the rank of the controllability matrix $Q = \big[B\ AB\ A^2B\ \cdots \ A^{n-1}B\big]=n$. 
\end{corollary}
\subsection{Linearization of Differential Systems}
Consider the nonlinear system
\begin{equation}
	\begin{aligned}
		\dot{x}({t})&=f(x(t),u(t)),\label{nonlin}
	\end{aligned}
\end{equation}
where the state x(t)is an n-dimensional vector, controller $u(t)$ is m-dimensional vector for all $t$, $f:\mathbb{R}^+\times \mathbb{R}^n \times \mathbb{R}^m \to \mathbb{R}^n$ is a non-linear function. 

Let $(x_0,u_0)$ be the reference point of the system \eqref{nonlin} then Taylor series expansion of the the function at the reference point is given by:
\begin{equation*}
	\begin{aligned}
	f(x_0+\delta x,u_0+\delta u)=f(x_0,u_0)+\frac{\partial f}{\partial x}\bigg|_{(x_0,u_0)} \delta x+\frac{\partial f}{\partial u}\bigg|_{(x_0,u_0)} \delta u+ higher\  order\ terms,
	\end{aligned}	
\end{equation*} 
and therefore we have:
\[\dot{x}_0+\delta \dot{x} \approx f(x_0,u_0)+\frac{\partial f}{\partial x}\bigg|_{(x_0,u_0)} \delta x+\frac{\partial f}{\partial u}\bigg|_{(x_0,u_0)}\delta u, \] 
simplifying, we get
\begin{equation}
	\delta x=\frac{\partial f}{\partial x}\bigg|_{(x_0,u_0)} \delta x+\frac{\partial f}{\partial u}\bigg|_{(x_0,u_0)}\delta u.\label{nonlin1}
\end{equation}
Define, $x=\delta x, u=\delta u, A=\frac{\partial f}{\partial x}\bigg|_{(x_0,u_0)}$ and $B=\frac{\partial f}{\partial u}\bigg|_{(x_0,u_0)}$ the system \eqref{nonlin1} becomes:
\begin{equation}
	\dot{x}=Ax+Bu. \label{nonlin3}
\end{equation}
The equation \eqref{nonlin3} is linear system corresponding to the system \eqref{nonlin}.

\section{Controllabiliy Analysis of the Motion of Satellite}
The equations of motion of satellite under the effect of oblateness of the earth is given by
\begin{equation}
 \ddot{\vec{r}}=-\frac{\mu}{r^3}\vec{r}+\vec{a}_O, \label{eqm1}
\end{equation}
where, $\mu=GM$, $G$ is gravitational constant and $M$ is mass of the earth and $\vec{a}_O$ is acceleration due to oblatenss of the earth, considering zonal harmonic $J_2$. The equations of motion in cylindrical coordinate systems represented by Humi\cite{Humi},
\begin{eqnarray}
  \nonumber \ddot{r}-r\dot{\theta}^2&=&-\mu r\bigg[\frac{1}{(r^2+z^2)^{\frac{3}{2}}}+\frac{3R^2J_2(r^2-4z^2)}{2(r^2+z^2)^{\frac{7}{2}}}\bigg],\\
  r\ddot{\theta}+2\dot{r}\dot{\theta}&=&0 \label{eqm2},\\
  \nonumber \ddot{z}&=&-\mu z\bigg[\frac{1}{(r^2+z^2)^{\frac{3}{2}}}+\frac{3R^2J_(3r^2-2z^2)}{2(r^2+z^2)^{\frac{7}{2}}}\bigg].
\end{eqnarray}
Under the effect of zonal  harmonic $J_{2}$, the satellite will deviate from its desired orbit, hence its motion becomes uncontrollable. Eventually it will hit on Earth. Hence, to control the motion of satellite we need to impose the controllers in the form of thrusters. Let $u_{1},\;u_{2}$ and $u_{3}$ represents thrusters in the $r,\;\theta$ and $z$ directions respectively. We analysed seven different cases viz. applying thurster(s) in 
\begin{enumerate}
	\item only $r$ direction,
	\item only $\theta$ direction,
	\item only $z$ direction,
	\item $r$ and $\theta$ direction,
	\item $r$ and $z$ direction,
	\item $\theta$ and $z$ direction,
	\item $r$, $\theta$ and $z$ direction.
\end{enumerate}
and check the controllability of system in each case.\\\\
Further we assume that the orbit of the satellite is circular with reference radius $\sigma$ and the angle $\theta=\omega t$. Since we have well established theory of controllability for first order system, we apply the following transformation to the system (\ref{eqm2}) after adding controllers in various directions to reduce it to a system of first order equations,
\begin{eqnarray}\label{Trans1}
	X_{1} & = & r-\sigma,\nonumber\\
	X_{2} & = & \dot{r},\nonumber\\
	X_{3} & = & \sigma\left(\theta-\omega t \right),\\
	X_{4} & = & \sigma\left(\dot{\theta}-\omega \right),\nonumber\\
	X_{5} & = & z,\nonumber\\
	X_{6} & = & \dot{z}.\nonumber
\end{eqnarray}
The study of controllability after applying thrusters in the different directions are discussed below.\\\\
\subsection{Adding the thruster $u_1(t)$ only in $r$ direction, the system (\ref{eqm2}) becomes:}
\begin{eqnarray}\label{eqm3}
	\ddot{r}-r\dot{\theta^{2}} & = & -\mu r\left[\frac{1}{\left(r^2+z^2 \right)^{3/2}}+\frac{3R^{2}J_{2}\left(r^{2}-4z^{2} \right)}{2\left(r^{2}+z^{2} \right)^{7/2}} \right]+u_{1}(t),\nonumber\\
	r\ddot{\theta}+2\dot{r}\dot{\theta} & = & 0,\\
	\ddot{z} & = & -\mu z\left[\frac{1}{\left(r^{2}+z^{2} \right)^{3/2}}+\frac{3R^{2}J_{2}\left(3r^{2}-2z^{2} \right)}{2\left(r^{2}+z^{2} \right)^{7/2}} \right].\nonumber
\end{eqnarray}
By transformation (\ref{Trans1}), system (\ref{eqm3}) takes the form:
\begin{eqnarray}\label{eqm4}
	\frac{dX_{1}}{dt} & = & X_{2},\nonumber\\
	\frac{dX_{2}}{dt} & = & \left(X_{1}+\sigma \right)\left(\frac{X_{4}}{\sigma}+\omega \right)^{2}-\mu\left(X_{1}+\sigma \right)\left\{\frac{1}{\left[\left(X_{1}+\sigma \right)^{2}+X_{5}^2 \right]^{3/2}}+\frac{3R^{2}J_{2}\left[\left(X_{1}+\sigma \right)^{2}-4X_{5}^{2} \right]}{2\left[\left(X_{1}+\sigma \right)^2+X_{5} \right]^{7/2}} \right\}+u_{1}\left(t \right),\nonumber\\
	\frac{dX_{3}}{dt} & = & X_{4},\\
	\frac{dX_{4}}{dt} & = & -\frac{2X_{5}\sigma\left(\frac{X_{4}}{\sigma}+\omega \right)}{\left(X_{1}+\sigma \right)},\nonumber\\
	\frac{dX_{5}}{dt} & = & X_{6},\nonumber\\
	\frac{dX_{6}}{dt} & = & -\mu X_{5}\left\{\frac{1}{\left[\left(X_{1}+\sigma \right)^{2}+X_{5}^2 \right]^{3/2}}+\frac{3R^{2}J_{2}\left[3\left(X_{1}+\sigma \right)^{2}-2X_{5}^{2} \right]}{2\left[\left(X_{1}+\sigma \right)^2+X_{5} \right]^{7/2}} \right\}.\nonumber
\end{eqnarray}
Now we linearize the system (\ref{eqm4}) about origin, we take
\begin{eqnarray}\label{f1}
	f_{1} & = &X_{2},\nonumber\\
	f_{2} & = &\left(X_{1}+\sigma \right)\left(\frac{X_{4}}{\sigma}+\omega \right)^{2}-\mu\left(X_{1}+\sigma \right)\left\{\frac{1}{\left[\left(X_{1}+\sigma \right)^{2}+X_{5}^2 \right]^{3/2}}+\frac{3R^{2}J_{2}\left[\left(X_{1}+\sigma \right)^{2}-4X_{5}^{2} \right]}{2\left[\left(X_{1}+\sigma \right)^2+X_{5} \right]^{7/2}} \right\}+u_{1}\left(t \right),\nonumber\\
	f_{3} & = & X_{4},\nonumber\\
	f_{4} & = & -\frac{2X_{5}\sigma\left(\frac{X_{4}}{\sigma}+\omega \right)}{\left(X_{1}+\sigma \right)},\nonumber\\
	f_{5} & = & X_{6},\nonumber\\
	f_{6} & = & -\mu X_{5}\left\{\frac{1}{\left[\left(X_{1}+\sigma \right)^{2}+X_{5}^2 \right]^{3/2}}+\frac{3R^{2}J_{2}\left[3\left(X_{1}+\sigma \right)^{2}-2X_{5}^{2} \right]}{2\left[\left(X_{1}+\sigma \right)^2+X_{5} \right]^{7/2}} \right\},\nonumber
\end{eqnarray}
therefore system (\ref{eqm4}) takes the form 
\[\dot{X}=AX+BU,\] 
where, $\dot{X}=\begin{bmatrix} \frac{dX_{1}}{dt} & \frac{dX_{2}}{dt} & \frac{dX_{3}}{dt} & \frac{dX_{4}}{dt} & \frac{dX_{5}}{dt} & \frac{dX_{6}}{dt}\end{bmatrix}',$ $ A=\begin{bmatrix} \frac{\partial\left(f_{1},\;f_{2},\;f_{3},\;f_{4},\;f_{5},\;f_{6} \right)}{\partial\left(X_{1},\;X_{2},\;X_{3},\;X_{4},\;X_{5},\;X_{6}\right)}\end{bmatrix}$ 
at origin,
$X=\begin{bmatrix}X_{1}\;X_{2}\;X_{3}\;X_{4}\;X_{5}\;X_{6}\end{bmatrix}'$,
$B=\begin{bmatrix}\frac{\partial f_{1}}{\partial u_{1}}\;\frac{\partial f_{2}}{\partial u_{1}}\;\frac{\partial f_{3}}{\partial u_{1}}\;\frac{\partial f_{4}}{\partial u_{1}}\;\frac{\partial f_{5}}{\partial u_{1}}\;\frac{\partial f_{6}}{\partial u_{1}} \end{bmatrix}'$ at origin and $u=\begin{bmatrix}u_{1}\end{bmatrix}$.
The values of $A$ and $B$ are
\begin{equation}\label{A}
A=\begin{bmatrix} 
0 & 1 & 0 & 0 & 0 & 0\\
1.000002542612694 & 0 & 0 & 2 & 0 & 0\\
0 & 0 & 0 & 1 & 0 & 0\\
0 & -0.000294117647059 & 0 & 0 & 0 & 0\\
0 & 0 & 0 & 0 & 0 & 1\\
0 & 0 & 0 & 0 & -0.00000127311747 & 0
\end{bmatrix},
\end{equation}
and $B=\begin{bmatrix}0\;1\;0\;0\;0\;0 \end{bmatrix}'$. The controllability matrix Q is given by
\begin{equation*}
 Q=\begin{bmatrix}B\;AB\;A^{2}B\;A^{3}B\;A^{4}B\;A^{5}B \end{bmatrix}=\begin{bmatrix}0 & 1 & 0 & -3 & 0 & 9\\1 & 0 & -3 & 0 & 9 & 0\\0 & 0 & -2 & 0 & 6 & 0\\0 & -2 & 0 & 6 & 0 & -18\\0 & 0 & 0 & 0 & 0 & 0\\0 & 0 & 0 & 0 & 0 & 0 \end{bmatrix}.
\end{equation*}
The rank of the matrix $Q$ is 3, which is not equal to the dimensions of the state $X\left(=6\right)$. By the Kalman's condition, the system is not controllable if we add the thruster only in radial direction $r$.
\subsection{Adding the thruster $u_1(t)$ only in $\theta$ direction, the system (\ref{eqm2}) becomes:}
\begin{eqnarray}\label{eqm5}
	\ddot{r}-r\dot{\theta^{2}} & = & -\mu r\left[\frac{1}{\left(r^2+z^2 \right)^{3/2}}+\frac{3R^{2}J_{2}\left(r^{2}-4z^{2} \right)}{2\left(r^{2}+z^{2} \right)^{7/2}} \right],\nonumber\\
	r\ddot{\theta}+2\dot{r}\dot{\theta} & = & u_{1}(t),\\
	\ddot{z} & = & -\mu z\left[\frac{1}{\left(r^{2}+z^{2} \right)^{3/2}}+\frac{3R^{2}J_{2}\left(3r^{2}-2z^{2} \right)}{2\left(r^{2}+z^{2} \right)^{7/2}} \right].\nonumber
\end{eqnarray}
By transformation (\ref{Trans1}), system (\ref{eqm5}) takes the form:
\begin{eqnarray}\label{eqm6}
	\frac{dX_{1}}{dt} & = & X_{2},\nonumber\\
	\frac{dX_{2}}{dt} & = & \left(X_{1}+\sigma \right)\left(\frac{X_{4}}{\sigma}+\omega \right)^{2}-\mu\left(X_{1}+\sigma \right)\left\{\frac{1}{\left[\left(X_{1}+\sigma \right)^{2}+X_{5}^2 \right]^{3/2}}+\frac{3R^{2}J_{2}\left[\left(X_{1}+\sigma \right)^{2}-4X_{5}^{2} \right]}{2\left[\left(X_{1}+\sigma \right)^2+X_{5} \right]^{7/2}} \right\},\nonumber\\
	\frac{dX_{3}}{dt} & = & X_{4},\\
	\frac{dX_{4}}{dt} & = & -\frac{2X_{5}\sigma\left(\frac{X_{4}}{\sigma}+\omega \right)}{\left(X_{1}+\sigma \right)}+u_{1}\left(t \right),\nonumber\\
	\frac{dX_{5}}{dt} & = & X_{6},\nonumber\\
	\frac{dX_{6}}{dt} & = & -\mu X_{5}\left\{\frac{1}{\left[\left(X_{1}+\sigma \right)^{2}+X_{5}^2 \right]^{3/2}}+\frac{3R^{2}J_{2}\left[3\left(X_{1}+\sigma \right)^{2}-2X_{5}^{2} \right]}{2\left[\left(X_{1}+\sigma \right)^2+X_{5} \right]^{7/2}} \right\}.\nonumber
\end{eqnarray}
We linearize the system (\ref{eqm6}) about origin by taking
\begin{eqnarray}\label{f2}
	f_{1} & = &X_{2},\nonumber\\
	f_{2} & = &\left(X_{1}+\sigma \right)\left(\frac{X_{4}}{\sigma}+\omega \right)^{2}-\mu\left(X_{1}+\sigma \right)\left\{\frac{1}{\left[\left(X_{1}+\sigma \right)^{2}+X_{5}^2 \right]^{3/2}}+\frac{3R^{2}J_{2}\left[\left(X_{1}+\sigma \right)^{2}-4X_{5}^{2} \right]}{2\left[\left(X_{1}+\sigma \right)^2+X_{5} \right]^{7/2}} \right\},\nonumber\\
	f_{3} & = & X_{4},\nonumber\\
	f_{4} & = & -\frac{2X_{5}\sigma\left(\frac{X_{4}}{\sigma}+\omega \right)}{\left(X_{1}+\sigma \right)}+u_{1}\left(t \right),\nonumber\\
	f_{5} & = & X_{6},\nonumber\\
	f_{6} & = & -\mu X_{5}\left\{\frac{1}{\left[\left(X_{1}+\sigma \right)^{2}+X_{5}^2 \right]^{3/2}}+\frac{3R^{2}J_{2}\left[3\left(X_{1}+\sigma \right)^{2}-2X_{5}^{2} \right]}{2\left[\left(X_{1}+\sigma \right)^2+X_{5} \right]^{7/2}} \right\},\nonumber
\end{eqnarray}
therefore the system (\ref{eqm6}) takes the form \[\dot{X}=AX+BU,\] where, $\dot{X}=\begin{bmatrix} \frac{dX_{1}}{dt} & \frac{dX_{2}}{dt} & \frac{dX_{3}}{dt} & \frac{dX_{4}}{dt} & \frac{dX_{5}}{dt} & \frac{dX_{6}}{dt}\end{bmatrix}'$, $A=\begin{bmatrix}\frac{\partial\left(f_{1},\;f_{2},\;f_{3},\;f_{4},\;f_{5},\;f_{6} \right)}{\partial\left(X_{1},\;X_{2},\;X_{3},\;X_{4},\;X_{5},\;X_{6}\right)}\end{bmatrix}$, at origin, 
$ X=\begin{bmatrix}X_{1}\;X_{2}\;X_{3}\;X_{4}\;X_{5}\;X_{6}\end{bmatrix}'$, $B=\begin{bmatrix}\frac{\partial f_{1}}{\partial u_{1}}\;\frac{\partial f_{2}}{\partial u_{1}}\;\frac{\partial f_{3}}{\partial u_{1}}\;\frac{\partial f_{4}}{\partial u_{1}}\;\frac{\partial f_{5}}{\partial u_{1}}\;\frac{\partial f_{6}}{\partial u_{1}} \end{bmatrix}'$ at origin, and  $u=\begin{bmatrix}u_{1}\end{bmatrix}.$ The matrix $A$ is given by (\ref{A}) and $B=\begin{bmatrix}0\;0\;0\;1\;0\;0 \end{bmatrix}'$. The controllability matrix $Q$ is given by  
\begin{equation*}
 Q=\begin{bmatrix}B\;AB\;A^{2}B\;A^{3}B\;A^{4}B\;A^{5}B \end{bmatrix}=\begin{bmatrix}0 & 0 & 2 & 0 & -6 & 0\\0 & 2 & 0 & -6 & 0 & 18\\0 & 1 & 0 & -4 & 0 & 12\\1 & 0 & -4 & 0 & 12 & 0\\0 & 0 & 0 & 0 & 0 & 0\\0 & 0 & 0 & 0 & 0 & 0 \end{bmatrix},
\end{equation*}
and the rank of $Q$ is 4, which is not equal to the dimensions of the state $X\left(=6\right)$. By the Kalman's condition, the system is not controllable if we add the thruster only in $\theta$ direction.
\subsection{Adding the thruster $u_{1}(t)$ only in $z$ direction}
The system (\ref{eqm2}) is written as
\begin{eqnarray}\label{eqm7}
	\ddot{r}-r\dot{\theta^{2}} & = & -\mu r\left[\frac{1}{\left(r^2+z^2 \right)^{3/2}}+\frac{3R^{2}J_{2}\left(r^{2}-4z^{2} \right)}{2\left(r^{2}+z^{2} \right)^{7/2}} \right],\nonumber\\
	r\ddot{\theta}+2\dot{r}\dot{\theta} & = 0,\\
	\ddot{z} & = & -\mu z\left[\frac{1}{\left(r^{2}+z^{2} \right)^{3/2}}+\frac{3R^{2}J_{2}\left(3r^{2}-2z^{2} \right)}{2\left(r^{2}+z^{2} \right)^{7/2}} \right]+u_{1}(t).\nonumber
\end{eqnarray}
By transformation (\ref{Trans1}), system (\ref{eqm7}) takes the form:
\begin{eqnarray}\label{eqm8}
	\frac{dX_{1}}{dt} & = & X_{2},\nonumber\\
	\frac{dX_{2}}{dt} & = & \left(X_{1}+\sigma \right)\left(\frac{X_{4}}{\sigma}+\omega \right)^{2}-\mu\left(X_{1}+\sigma \right)\left\{\frac{1}{\left[\left(X_{1}+\sigma \right)^{2}+X_{5}^2 \right]^{3/2}}+\frac{3R^{2}J_{2}\left[\left(X_{1}+\sigma \right)^{2}-4X_{5}^{2} \right]}{2\left[\left(X_{1}+\sigma \right)^2+X_{5} \right]^{7/2}} \right\},\nonumber\\
	\frac{dX_{3}}{dt} & = & X_{4},\\
	\frac{dX_{4}}{dt} & = & -\frac{2X_{5}\sigma\left(\frac{X_{4}}{\sigma}+\omega \right)}{\left(X_{1}+\sigma \right)},\nonumber\\
	\frac{dX_{5}}{dt} & = & X_{6},\nonumber\\
	\frac{dX_{6}}{dt} & = & -\mu X_{5}\left\{\frac{1}{\left[\left(X_{1}+\sigma \right)^{2}+X_{5}^2 \right]^{3/2}}+\frac{3R^{2}J_{2}\left[3\left(X_{1}+\sigma \right)^{2}-2X_{5}^{2} \right]}{2\left[\left(X_{1}+\sigma \right)^2+X_{5} \right]^{7/2}} \right\}+u_{1}\left(t \right).\nonumber
\end{eqnarray}
Linearizing the system (\ref{eqm8}) about origin by taking
\begin{eqnarray}\label{f3}
	f_{1} & = &X_{2},\nonumber\\
	f_{2} & = &\left(X_{1}+\sigma \right)\left(\frac{X_{4}}{\sigma}+\omega \right)^{2}-\mu\left(X_{1}+\sigma \right)\left\{\frac{1}{\left[\left(X_{1}+\sigma \right)^{2}+X_{5}^2 \right]^{3/2}}+\frac{3R^{2}J_{2}\left[\left(X_{1}+\sigma \right)^{2}-4X_{5}^{2} \right]}{2\left[\left(X_{1}+\sigma \right)^2+X_{5} \right]^{7/2}} \right\},\nonumber\\
	f_{3} & = & X_{4},\nonumber\\
	f_{4} & = & -\frac{2X_{5}\sigma\left(\frac{X_{4}}{\sigma}+\omega \right)}{\left(X_{1}+\sigma \right)},\nonumber\\
	f_{5} & = & X_{6},\nonumber\\
	f_{6} & = & -\mu X_{5}\left\{\frac{1}{\left[\left(X_{1}+\sigma \right)^{2}+X_{5}^2 \right]^{3/2}}+\frac{3R^{2}J_{2}\left[3\left(X_{1}+\sigma \right)^{2}-2X_{5}^{2} \right]}{2\left[\left(X_{1}+\sigma \right)^2+X_{5} \right]^{7/2}} \right\}+u_{1}\left(t \right),\nonumber
\end{eqnarray}
and the system (\ref{eqm8}) takes the form \[\dot{X}=AX+BU,\] where,
$\dot{X}=\begin{bmatrix} \frac{dX_{1}}{dt} & \frac{dX_{2}}{dt} & \frac{dX_{3}}{dt} & \frac{dX_{4}}{dt} & \frac{dX_{5}}{dt} & \frac{dX_{6}}{dt}\end{bmatrix}'$, $A=\begin{bmatrix}\frac{\partial\left(f_{1},\;f_{2},\;f_{3},\;f_{4},\;f_{5},\;f_{6} \right)}{\partial\left(X_{1},\;X_{2},\;X_{3},\;X_{4},\;X_{5},\;X_{6}\right)}\end{bmatrix}$ at origin, $\begin{bmatrix}X=X_{1}\;X_{2}\;X_{3}\;X_{4}\;X_{5}\;X_{6}\end{bmatrix}'$, $B=\begin{bmatrix}\frac{\partial f_{1}}{\partial u_{1}}\;\frac{\partial f_{2}}{\partial u_{1}}\;\frac{\partial f_{3}}{\partial u_{1}}\;\frac{\partial f_{4}}{\partial u_{1}}\;\frac{\partial f_{5}}{\partial u_{1}}\;\frac{\partial f_{6}}{\partial u_{1}} \end{bmatrix}'$ at origin and $U=\begin{bmatrix}u_{1} \end{bmatrix}$. The values of $A$ as (\ref{A}) and $B=\begin{bmatrix}0\;0\;0\;0\;0\;1 \end{bmatrix}'$. The controllability matrix $Q$ is given by 
\begin{equation*}
 Q=\begin{bmatrix}B\;AB\;A^{2}B\;A^{3}B\;A^{4}B\;A^{5}B \end{bmatrix}=\begin{bmatrix}0 & 0 & 0 & 0 & 0 & 0\\0 & 0 & 0 & 0 & 0 & 0\\0 & 0 & 0 & 0 & 0 & 0\\0 & 0 & 0 & 0 & 0 & 0\\0 & 1 & 0 & 0 & 0 & 0\\1 & 0 & 0 & 0 & 0 & 0 \end{bmatrix},
\end{equation*} 
Therefore, rank of the matrix $Q=2$ which is not equal to the dimensions of the state $X\left(=6\right)$. By the Kalman's condition, the system is not controllable if we add the thruster only in $z$ direction.
\subsection{Adding thrusters $u_1(t)$  and $u_2(t)$ in $r$ and $\theta$ direction:}
The system (\ref{eqm2}) becomes:
\begin{eqnarray}\label{eqm9}
	\ddot{r}-r\dot{\theta^{2}} & = & -\mu r\left[\frac{1}{\left(r^2+z^2 \right)^{3/2}}+\frac{3R^{2}J_{2}\left(r^{2}-4z^{2} \right)}{2\left(r^{2}+z^{2} \right)^{7/2}} \right]+u_{1}(t),\nonumber\\
	r\ddot{\theta}+2\dot{r}\dot{\theta} & = u_{2}(t),\\
	\ddot{z} & = & -\mu z\left[\frac{1}{\left(r^{2}+z^{2} \right)^{3/2}}+\frac{3R^{2}J_{2}\left(3r^{2}-2z^{2} \right)}{2\left(r^{2}+z^{2} \right)^{7/2}} \right].\nonumber
\end{eqnarray}
By transformation (\ref{Trans1}), system (\ref{eqm9}) takes the form
\begin{eqnarray}\label{eqm10}
	\frac{dX_{1}}{dt} & = & X_{2},\nonumber\\
	\frac{dX_{2}}{dt} & = & \left(X_{1}+\sigma \right)\left(\frac{X_{4}}{\sigma}+\omega \right)^{2}-\mu\left(X_{1}+\sigma \right)\left\{\frac{1}{\left[\left(X_{1}+\sigma \right)^{2}+X_{5}^2 \right]^{3/2}}+\frac{3R^{2}J_{2}\left[\left(X_{1}+\sigma \right)^{2}-4X_{5}^{2} \right]}{2\left[\left(X_{1}+\sigma \right)^2+X_{5} \right]^{7/2}} \right\}+u_{1}(t),\nonumber\\
	\frac{dX_{3}}{dt} & = & X_{4},\\
	\frac{dX_{4}}{dt} & = & -\frac{2X_{5}\sigma\left(\frac{X_{4}}{\sigma}+\omega \right)}{\left(X_{1}+\sigma \right)}+u_{2}(t),\nonumber\\
	\frac{dX_{5}}{dt} & = & X_{6},\nonumber\\
	\frac{dX_{6}}{dt} & = & -\mu X_{5}\left\{\frac{1}{\left[\left(X_{1}+\sigma \right)^{2}+X_{5}^2 \right]^{3/2}}+\frac{3R^{2}J_{2}\left[3\left(X_{1}+\sigma \right)^{2}-2X_{5}^{2} \right]}{2\left[\left(X_{1}+\sigma \right)^2+X_{5} \right]^{7/2}} \right\}.\nonumber
\end{eqnarray}
For linearizing the system (\ref{eqm10}) about origin, we take
\begin{eqnarray}\label{f4}
	f_{1} & = &X_{2},\nonumber\\
	f_{2} & = &\left(X_{1}+\sigma \right)\left(\frac{X_{4}}{\sigma}+\omega \right)^{2}-\mu\left(X_{1}+\sigma \right)\left\{\frac{1}{\left[\left(X_{1}+\sigma \right)^{2}+X_{5}^2 \right]^{3/2}}+\frac{3R^{2}J_{2}\left[\left(X_{1}+\sigma \right)^{2}-4X_{5}^{2} \right]}{2\left[\left(X_{1}+\sigma \right)^2+X_{5} \right]^{7/2}} \right\}+u_{1}(t),\nonumber\\
	f_{3} & = & X_{4},\nonumber\\
	f_{4} & = & -\frac{2X_{5}\sigma\left(\frac{X_{4}}{\sigma}+\omega \right)}{\left(X_{1}+\sigma \right)}+u_{2}(t),\nonumber\\
	f_{5} & = & X_{6},\nonumber\\
	f_{6} & = & -\mu X_{5}\left\{\frac{1}{\left[\left(X_{1}+\sigma \right)^{2}+X_{5}^2 \right]^{3/2}}+\frac{3R^{2}J_{2}\left[3\left(X_{1}+\sigma \right)^{2}-2X_{5}^{2} \right]}{2\left[\left(X_{1}+\sigma \right)^2+X_{5} \right]^{7/2}} \right\},\nonumber
\end{eqnarray}
Therefore the system (\ref{eqm10}) take the form \[\dot{X}=AX+BU,\] where, 
$\dot{X}=\begin{bmatrix} \frac{dX_{1}}{dt} & \frac{dX_{2}}{dt} & \frac{dX_{3}}{dt} & \frac{dX_{4}}{dt} & \frac{dX_{5}}{dt} & \frac{dX_{6}}{dt}\end{bmatrix}'$, $A=\begin{bmatrix}\frac{\partial\left(f_{1},\;f_{2},\;f_{3},\;f_{4},\;f_{5},\;f_{6} \right)}{\partial\left(X_{1},\;X_{2},\;X_{3},\;X_{4},\;X_{5},\;X_{6}\right)}\end{bmatrix}$ at origin, $X=\begin{bmatrix}X_{1}\;X_{2}\;X_{3}\;X_{4}\;X_{5}\;X_{6}\end{bmatrix}'$, $B=\begin{bmatrix}\frac{\partial f_{1}}{\partial u_{1}}\;\frac{\partial f_{2}}{\partial u_{1}}\;\frac{\partial f_{3}}{\partial u_{1}}\;\frac{\partial f_{4}}{\partial u_{1}}\;\frac{\partial f_{5}}{\partial u_{1}}\;\frac{\partial f_{6}}{\partial u_{1}}\\
\frac{\partial f_{1}}{\partial u_{2}}\;\frac{\partial f_{2}}{\partial u_{2}}\;\frac{\partial f_{3}}{\partial u_{2}}\;\frac{\partial f_{4}}{\partial u_{2}}\;\frac{\partial f_{5}}{\partial u_{2}}\;\frac{\partial f_{6}}{\partial u_{2}}\end{bmatrix}'$ at origin and $u=\begin{bmatrix}u_{1}\; u_{2} \end{bmatrix}'$. We obtain the values of $A$ as (\ref{A})and $B=\begin{bmatrix}0\;1\;0\;0\;0\;0\\0\;0\;0\;1\;0\;0 \end{bmatrix}'.$ The controllability matrix $Q$ is given by
\setcounter{MaxMatrixCols}{20}
\begin{equation*}
 Q=\begin{bmatrix}B\;AB\;A^{2}B\;A^{3}B\;A^{4}B\;A^{5}B \end{bmatrix}=\begin{bmatrix}0 & 0 & 1 & 0 & 0 & 2 & -3 & 0 & 0 & -6 & 9 & 0\\
1 & 0 & 0 & 2 & -3 & 0 & 0 & -6 & 9 & 0 & 0 & 18\\
0 & 0 & 0 & 1 & -2 & 0 & 0 & -4 & 6 & 0 & 0 & 12\\
0 & 1 & -2 & 0 & 0 & -4 & 6 & 0 & 0 & 12 & -18 & 0\\
0 & 0 & 0 & 0 & 0 & 0 & 0 & 0 & 0 & 0 & 0 & 0\\
0 & 0 & 0 & 0 & 0 & 0 & 0 & 0 & 0 & 0 & 0 & 0
\end{bmatrix},
\end{equation*}

The rank of the matrix $Q=4$ which is not equal to the dimensions of the state $X\left(=6\right)$. By the Kalman's condition, we conclude that the system is not controlllable if we add the thrusters in $r$ and $\theta$ direction.
\subsection{Adding thrusters $u_1(t)$ and $u_2(t)$ in $r$ and $z$ direction:}
The system (\ref{eqm2}) becomes:
\begin{eqnarray}\label{eqm11}
	\ddot{r}-r\dot{\theta^{2}} & = & -\mu r\left[\frac{1}{\left(r^2+z^2 \right)^{3/2}}+\frac{3R^{2}J_{2}\left(r^{2}-4z^{2} \right)}{2\left(r^{2}+z^{2} \right)^{7/2}} \right]+u_{1}(t),\nonumber\\
	r\ddot{\theta}+2\dot{r}\dot{\theta} & = 0,\\
	\ddot{z} & = & -\mu z\left[\frac{1}{\left(r^{2}+z^{2} \right)^{3/2}}+\frac{3R^{2}J_{2}\left(3r^{2}-2z^{2} \right)}{2\left(r^{2}+z^{2} \right)^{7/2}} \right]+u_{2}(t).\nonumber
\end{eqnarray}
By transformation (\ref{Trans1}), system (\ref{eqm11}) takes the form
\begin{eqnarray}\label{eqm12}
	\frac{dX_{1}}{dt} & = & X_{2},\nonumber\\
	\frac{dX_{2}}{dt} & = & \left(X_{1}+\sigma \right)\left(\frac{X_{4}}{\sigma}+\omega \right)^{2}-\mu\left(X_{1}+\sigma \right)\left\{\frac{1}{\left[\left(X_{1}+\sigma \right)^{2}+X_{5}^2 \right]^{3/2}}+\frac{3R^{2}J_{2}\left[\left(X_{1}+\sigma \right)^{2}-4X_{5}^{2} \right]}{2\left[\left(X_{1}+\sigma \right)^2+X_{5} \right]^{7/2}} \right\}+u_{1}(t),\nonumber\\
	\frac{dX_{3}}{dt} & = & X_{4},\\
	\frac{dX_{4}}{dt} & = & -\frac{2X_{5}\sigma\left(\frac{X_{4}}{\sigma}+\omega \right)}{\left(X_{1}+\sigma \right)},\nonumber\\
	\frac{dX_{5}}{dt} & = & X_{6},\nonumber\\
	\frac{dX_{6}}{dt} & = & -\mu X_{5}\left\{\frac{1}{\left[\left(X_{1}+\sigma \right)^{2}+X_{5}^2 \right]^{3/2}}+\frac{3R^{2}J_{2}\left[3\left(X_{1}+\sigma \right)^{2}-2X_{5}^{2} \right]}{2\left[\left(X_{1}+\sigma \right)^2+X_{5} \right]^{7/2}} \right\}+u_{2}(t).\nonumber
\end{eqnarray}
For linearizing the system (\ref{eqm12}) about origin, we take
\begin{eqnarray}\label{f5}
	f_{1} & = &X_{2},\nonumber\\
	f_{2} & = &\left(X_{1}+\sigma \right)\left(\frac{X_{4}}{\sigma}+\omega \right)^{2}-\mu\left(X_{1}+\sigma \right)\left\{\frac{1}{\left[\left(X_{1}+\sigma \right)^{2}+X_{5}^2 \right]^{3/2}}+\frac{3R^{2}J_{2}\left[\left(X_{1}+\sigma \right)^{2}-4X_{5}^{2} \right]}{2\left[\left(X_{1}+\sigma \right)^2+X_{5} \right]^{7/2}} \right\}+u_{1}(t),\nonumber\\
	f_{3} & = & X_{4},\nonumber\\
	f_{4} & = & -\frac{2X_{5}\sigma\left(\frac{X_{4}}{\sigma}+\omega \right)}{\left(X_{1}+\sigma \right)},\nonumber\\
	f_{5} & = & X_{6},\nonumber\\
	f_{6} & = & -\mu X_{5}\left\{\frac{1}{\left[\left(X_{1}+\sigma \right)^{2}+X_{5}^2 \right]^{3/2}}+\frac{3R^{2}J_{2}\left[3\left(X_{1}+\sigma \right)^{2}-2X_{5}^{2} \right]}{2\left[\left(X_{1}+\sigma \right)^2+X_{5} \right]^{7/2}} \right\}+u_{2}(t),\nonumber
\end{eqnarray}
Therefore the system (\ref{eqm12}) takes the form \[\dot{X}=AX+BU,\] where, $\dot{X}=\begin{bmatrix} \frac{dX_{1}}{dt} & \frac{dX_{2}}{dt} & \frac{dX_{3}}{dt} & \frac{dX_{4}}{dt} & \frac{dX_{5}}{dt} & \frac{dX_{6}}{dt}\end{bmatrix}'$, $A=\begin{bmatrix}\frac{\partial\left(f_{1},\;f_{2},\;f_{3},\;f_{4},\;f_{5},\;f_{6} \right)}{\partial\left(X_{1},\;X_{2},\;X_{3},\;X_{4},\;X_{5},\;X_{6}\right)}\end{bmatrix}$ at origin, $X=\begin{bmatrix}X_{1}\;X_{2}\;X_{3}\;X_{4}\;X_{5}\;X_{6}\end{bmatrix}'$, $B=\begin{bmatrix}\frac{\partial f_{1}}{\partial u_{1}}\;\frac{\partial f_{2}}{\partial u_{1}}\;\frac{\partial f_{3}}{\partial u_{1}}\;\frac{\partial f_{4}}{\partial u_{1}}\;\frac{\partial f_{5}}{\partial u_{1}}\;\frac{\partial f_{6}}{\partial u_{1}}\\
\frac{\partial f_{1}}{\partial u_{2}}\;\frac{\partial f_{2}}{\partial u_{2}}\;\frac{\partial f_{3}}{\partial u_{2}}\;\frac{\partial f_{4}}{\partial u_{2}}\;\frac{\partial f_{5}}{\partial u_{2}}\;\frac{\partial f_{6}}{\partial u_{2}}
 \end{bmatrix}'$ at origin and $u=\begin{bmatrix}u_{1}\; u_{2} \end{bmatrix}'$,
 The matrix $A$ is given in (\ref{A}) and $B=\begin{bmatrix}0\;1\;0\;0\;0\;0\\0\;0\;0\;0\;0\;1 \end{bmatrix}'$. 
 The controllability matrix $Q$ is given by
\begin{equation*}
 Q= \begin{bmatrix}B\;AB\;A^{2}B\;A^{3}B\;A^{4}B\;A^{5}B \end{bmatrix}=\begin{bmatrix}
0 & 0 & 1 & 0 & 0 & 0 & -3 & 0 & 0 & 0 & 9 & 0\\
1 & 0 & 0 & 0 & -3 & 0 & 0 & 0 & 9 & 0 & 0 & 0\\
0 & 0 & 0 & 0 & -2 & 0 & 0 & 0 & 6 & 0 & 0 & 0\\
0 & 0 & -2 & 0 & 0 & 0 & 6 & 0 & 0 & 0 & -18 & 0\\
0 & 0 & 0 & 1 & 0 & 0 & 0 & 0 & 0 & 0 & 0 & 0\\
0 & 1 & 0 & 0 & 0 & 0 & 0 & 0 & 0 & 0 & 0 & 0
\end{bmatrix},
\end{equation*}
The rank of the matrix $Q$ is 5, which is not equal to the dimensions of the state $X\left(=6\right)$. By the Kalman's condition, we conclude that the system is not controllable if we add the thrusters in $r$ and $z$ direction.
\subsection{Adding thrusters  $u_1(t)$  and $u_2(t)$ in $\theta$ and $z$ direction:}
The system (\ref{eqm2}) becomes:
\begin{eqnarray}\label{eqm13}
	\ddot{r}-r\dot{\theta^{2}} & = & -\mu r\left[\frac{1}{\left(r^2+z^2 \right)^{3/2}}+\frac{3R^{2}J_{2}\left(r^{2}-4z^{2} \right)}{2\left(r^{2}+z^{2} \right)^{7/2}} \right],\nonumber\\
	r\ddot{\theta}+2\dot{r}\dot{\theta} & = u_{1}(t),\\
	\ddot{z} & = & -\mu z\left[\frac{1}{\left(r^{2}+z^{2} \right)^{3/2}}+\frac{3R^{2}J_{2}\left(3r^{2}-2z^{2} \right)}{2\left(r^{2}+z^{2} \right)^{7/2}} \right]+u_{2}(t).\nonumber
\end{eqnarray}
By transformation (\ref{Trans1}), system (\ref{eqm13}) takes the form
\begin{eqnarray}\label{eqm14}
	\frac{dX_{1}}{dt} & = & X_{2},\nonumber\\
	\frac{dX_{2}}{dt} & = & \left(X_{1}+\sigma \right)\left(\frac{X_{4}}{\sigma}+\omega \right)^{2}-\mu\left(X_{1}+\sigma \right)\left\{\frac{1}{\left[\left(X_{1}+\sigma \right)^{2}+X_{5}^2 \right]^{3/2}}+\frac{3R^{2}J_{2}\left[\left(X_{1}+\sigma \right)^{2}-4X_{5}^{2} \right]}{2\left[\left(X_{1}+\sigma \right)^2+X_{5} \right]^{7/2}} \right\},\nonumber\\
	\frac{dX_{3}}{dt} & = & X_{4},\\
	\frac{dX_{4}}{dt} & = & -\frac{2X_{5}\sigma\left(\frac{X_{4}}{\sigma}+\omega \right)}{\left(X_{1}+\sigma \right)}+u_{1}(t),\nonumber\\
	\frac{dX_{5}}{dt} & = & X_{6},\nonumber\\
	\frac{dX_{6}}{dt} & = & -\mu X_{5}\left\{\frac{1}{\left[\left(X_{1}+\sigma \right)^{2}+X_{5}^2 \right]^{3/2}}+\frac{3R^{2}J_{2}\left[3\left(X_{1}+\sigma \right)^{2}-2X_{5}^{2} \right]}{2\left[\left(X_{1}+\sigma \right)^2+X_{5} \right]^{7/2}} \right\}+u_{2}(t).\nonumber
\end{eqnarray}
For linearizing the system (\ref{eqm14}) about origin, we take
\begin{eqnarray}\label{f6}
	f_{1} & = &X_{2},\nonumber\\
	f_{2} & = &\left(X_{1}+\sigma \right)\left(\frac{X_{4}}{\sigma}+\omega \right)^{2}-\mu\left(X_{1}+\sigma \right)\left\{\frac{1}{\left[\left(X_{1}+\sigma \right)^{2}+X_{5}^2 \right]^{3/2}}+\frac{3R^{2}J_{2}\left[\left(X_{1}+\sigma \right)^{2}-4X_{5}^{2} \right]}{2\left[\left(X_{1}+\sigma \right)^2+X_{5} \right]^{7/2}} \right\},\nonumber\\
	f_{3} & = & X_{4},\nonumber\\
	f_{4} & = & -\frac{2X_{5}\sigma\left(\frac{X_{4}}{\sigma}+\omega \right)}{\left(X_{1}+\sigma \right)}+u_{1}(t),\nonumber\\
	f_{5} & = & X_{6},\nonumber\\
	f_{6} & = & -\mu X_{5}\left\{\frac{1}{\left[\left(X_{1}+\sigma \right)^{2}+X_{5}^2 \right]^{3/2}}+\frac{3R^{2}J_{2}\left[3\left(X_{1}+\sigma \right)^{2}-2X_{5}^{2} \right]}{2\left[\left(X_{1}+\sigma \right)^2+X_{5} \right]^{7/2}} \right\}+u_{2}(t),\nonumber
\end{eqnarray}
Therefore the system (\ref{eqm14}) takes the form: \[\dot{X}=AX+BU,\] where, $\dot{X}=\begin{bmatrix} \frac{dX_{1}}{dt} & \frac{dX_{2}}{dt} & \frac{dX_{3}}{dt} & \frac{dX_{4}}{dt} & \frac{dX_{5}}{dt} & \frac{dX_{6}}{dt}\end{bmatrix}'$, $A=\begin{bmatrix}\frac{\partial\left(f_{1},\;f_{2},\;f_{3},\;f_{4},\;f_{5},\;f_{6} \right)}{\partial\left(X_{1},\;X_{2},\;X_{3},\;X_{4},\;X_{5},\;X_{6}\right)}\end{bmatrix}$ at origin, $X=\begin{bmatrix}X_{1}\;X_{2}\;X_{3}\;X_{4}\;X_{5}\;X_{6}\end{bmatrix}'$, $B=\begin{bmatrix}\frac{\partial f_{1}}{\partial u_{1}}\;\frac{\partial f_{2}}{\partial u_{1}}\;\frac{\partial f_{3}}{\partial u_{1}}\;\frac{\partial f_{4}}{\partial u_{1}}\;\frac{\partial f_{5}}{\partial u_{1}}\;\frac{\partial f_{6}}{\partial u_{1}}\\
\frac{\partial f_{1}}{\partial u_{2}}\;\frac{\partial f_{2}}{\partial u_{2}}\;\frac{\partial f_{3}}{\partial u_{2}}\;\frac{\partial f_{4}}{\partial u_{2}}\;\frac{\partial f_{5}}{\partial u_{2}}\;\frac{\partial f_{6}}{\partial u_{2}}
 \end{bmatrix}'$ at origin and $u=\begin{bmatrix}u_{1}\; u_{2} \end{bmatrix}'$. The matrix $A$ is given in (\ref{A}) and $B=\begin{bmatrix}0\;0\;0\;1\;0\;0\\0\;0\;0\;0\;0\;1 \end{bmatrix}'$. The controllability matrix $Q$ is given by
 \begin{equation*}
 Q=\begin{bmatrix}B\;AB\;A^{2}B\;A^{3}B\;A^{4}B\;A^{5}B \end{bmatrix}=\begin{bmatrix}
0 & 0 & 0 & 0 & 2 & 0 & 0 & 0 & -6 & 0 & 0 & 0\\
0 & 0 & 2 & 0 & 0 & 0 & -6 & 0 & 0 & 0 & 18 & 0\\
0 & 0 & 1 & 0 & 0 & 0 & -4 & 0 & 0 & 0 & 12 & 0\\
1 & 0 & 0 & 0 & -4 & 0 & 0 & 0 & 12 & 0 & 0 & 0\\
0 & 0 & 0 & 1 & 0 & 0 & 0 & 0 & 0 & 0 & 0 & 0\\
0 & 1 & 0 & 0 & 0 & 0 & 0 & 0 & 0 & 0 & 0 & 0
\end{bmatrix},
\end{equation*} 
The rank of $Q$ is 6, which is equal to the dimensions of the state $X\left(=6\right)$. Hence by the Kalman's condition, we conclude that the system is controllable if we add thrusters $u_1(t)$ and $u_2(t)$ in $\theta$ and $z$ directions. The figure-1 shows, that the system is steered from the initial point 
$[1\; 2\; 3\; 4\; 5\; 6]'$ to the final point $[6\; 5\; 4\; 3\; 2\; 1]'$ during the time interval $[0 , 10]$, by applying the controllers, i.e. thrusters  $u_1(t)$  and $u_2(t)$ in $\theta$ and $z$ direction. 
\begin{figure*}[ht]
\begin{tabular}{cc}
\includegraphics[width=9cm,height=9cm]{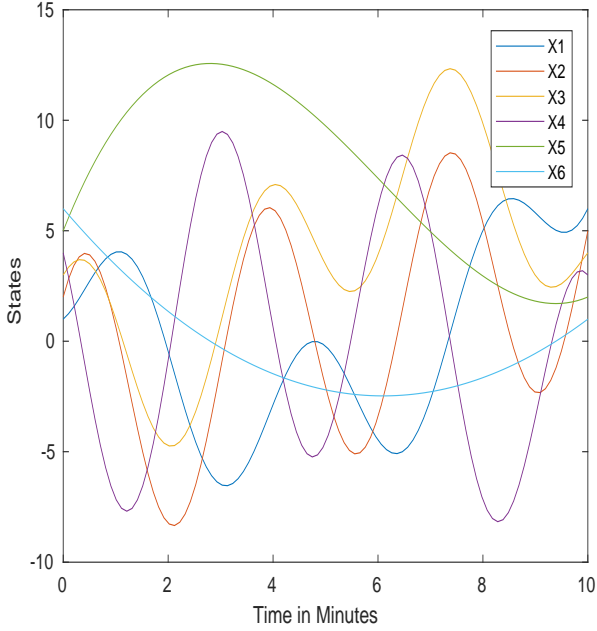} \end{tabular}
\caption{State Control of the System under the effect of zonal harmonic $J_2$}
\end{figure*}
The graph of the controllers i.e. thrusters  $u_1(t)$  and $u_2(t)$ in $\theta$ are shown in the figure-2:
\begin{figure*}[ht]
\begin{tabular}{cc}
\includegraphics[width=9cm,height=9cm]{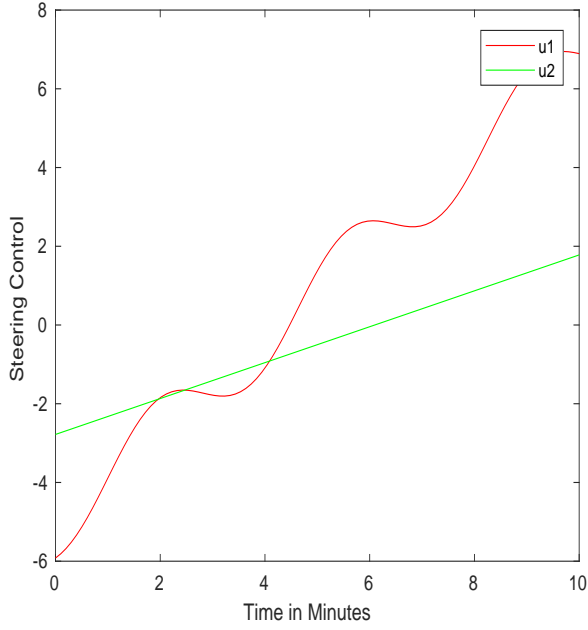} \end{tabular}
\caption{Steering Control of the System under the effect of zonal harmonic $J_2$}
\end{figure*}
\subsection{If we add the thrusters  in all the three directions i.e. $r$, $\theta$ and $z$ directions:}
The system (\ref{eqm2}) is written as
\begin{eqnarray}\label{eqm15}
	\ddot{r}-r\dot{\theta^{2}} & = & -\mu r\left[\frac{1}{\left(r^2+z^2 \right)^{3/2}}+\frac{3R^{2}J_{2}\left(r^{2}-4z^{2} \right)}{2\left(r^{2}+z^{2} \right)^{7/2}} \right]+u_{1}(t),\nonumber\\
	r\ddot{\theta}+2\dot{r}\dot{\theta} & = u_{2}(t),\\
	\ddot{z} & = & -\mu z\left[\frac{1}{\left(r^{2}+z^{2} \right)^{3/2}}+\frac{3R^{2}J_{2}\left(3r^{2}-2z^{2} \right)}{2\left(r^{2}+z^{2} \right)^{7/2}} \right]+u_{3}(t).\nonumber
\end{eqnarray}
By transformation (\ref{Trans1}), the system (\ref{eqm15}) takes the form
\begin{eqnarray}\label{eqm16}
	\frac{dX_{1}}{dt} & = & X_{2},\nonumber\\
	\frac{dX_{2}}{dt} & = & \left(X_{1}+\sigma \right)\left(\frac{X_{4}}{\sigma}+\omega \right)^{2}-\mu\left(X_{1}+\sigma \right)\left\{\frac{1}{\left[\left(X_{1}+\sigma \right)^{2}+X_{5}^2 \right]^{3/2}}+\frac{3R^{2}J_{2}\left[\left(X_{1}+\sigma \right)^{2}-4X_{5}^{2} \right]}{2\left[\left(X_{1}+\sigma \right)^2+X_{5} \right]^{7/2}} \right\}+u_{1}(t),\nonumber\\
	\frac{dX_{3}}{dt} & = & X_{4},\\
	\frac{dX_{4}}{dt} & = & -\frac{2X_{5}\sigma\left(\frac{X_{4}}{\sigma}+\omega \right)}{\left(X_{1}+\sigma \right)}+u_{2}(t),\nonumber\\
	\frac{dX_{5}}{dt} & = & X_{6},\nonumber\\
	\frac{dX_{6}}{dt} & = & -\mu X_{5}\left\{\frac{1}{\left[\left(X_{1}+\sigma \right)^{2}+X_{5}^2 \right]^{3/2}}+\frac{3R^{2}J_{2}\left[3\left(X_{1}+\sigma \right)^{2}-2X_{5}^{2} \right]}{2\left[\left(X_{1}+\sigma \right)^2+X_{5} \right]^{7/2}} \right\}+u_{3}(t).\nonumber
\end{eqnarray}
Now we linearize the system (\ref{eqm16}) about origin, we take
\begin{eqnarray}\label{f7}
	f_{1} & = &X_{2},\nonumber\\
	f_{2} & = &\left(X_{1}+\sigma \right)\left(\frac{X_{4}}{\sigma}+\omega \right)^{2}-\mu\left(X_{1}+\sigma \right)\left\{\frac{1}{\left[\left(X_{1}+\sigma \right)^{2}+X_{5}^2 \right]^{3/2}}+\frac{3R^{2}J_{2}\left[\left(X_{1}+\sigma \right)^{2}-4X_{5}^{2} \right]}{2\left[\left(X_{1}+\sigma \right)^2+X_{5} \right]^{7/2}} \right\}+u_{1}(t),\nonumber\\
	f_{3} & = & X_{4},\nonumber\\
	f_{4} & = & -\frac{2X_{5}\sigma\left(\frac{X_{4}}{\sigma}+\omega \right)}{\left(X_{1}+\sigma \right)}+u_{2}(t),\nonumber\\
	f_{5} & = & X_{6},\nonumber\\
	f_{6} & = & -\mu X_{5}\left\{\frac{1}{\left[\left(X_{1}+\sigma \right)^{2}+X_{5}^2 \right]^{3/2}}+\frac{3R^{2}J_{2}\left[3\left(X_{1}+\sigma \right)^{2}-2X_{5}^{2} \right]}{2\left[\left(X_{1}+\sigma \right)^2+X_{5} \right]^{7/2}} \right\}+u_{3}(t),\nonumber
\end{eqnarray}
and write the system (\ref{eqm16}) in the form \[\dot{X}=AX+BU,\] where $\dot{X}=\begin{bmatrix} \frac{dX_{1}}{dt} & \frac{dX_{2}}{dt} & \frac{dX_{3}}{dt} & \frac{dX_{4}}{dt} & \frac{dX_{5}}{dt} & \frac{dX_{6}}{dt}\end{bmatrix}'$, $A=\begin{bmatrix}\frac{\partial\left(f_{1},\;f_{2},\;f_{3},\;f_{4},\;f_{5},\;f_{6} \right)}{\partial\left(X_{1},\;X_{2},\;X_{3},\;X_{4},\;X_{5},\;X_{6}\right)}\end{bmatrix}$, at origin, $X=\begin{bmatrix}X_{1}\;X_{2}\;X_{3}\;X_{4}\;X_{5}\;X_{6}\end{bmatrix}'$, $B=\begin{bmatrix}\frac{\partial f_{1}}{\partial u_{1}}\;\frac{\partial f_{2}}{\partial u_{1}}\;\frac{\partial f_{3}}{\partial u_{1}}\;\frac{\partial f_{4}}{\partial u_{1}}\;\frac{\partial f_{5}}{\partial u_{1}}\;\frac{\partial f_{6}}{\partial u_{1}}\\
\frac{\partial f_{1}}{\partial u_{2}}\;\frac{\partial f_{2}}{\partial u_{2}}\;\frac{\partial f_{3}}{\partial u_{2}}\;\frac{\partial f_{4}}{\partial u_{2}}\;\frac{\partial f_{5}}{\partial u_{2}}\;\frac{\partial f_{6}}{\partial u_{2}}\\
\frac{\partial f_{1}}{\partial u_{3}}\;\frac{\partial f_{2}}{\partial u_{3}}\;\frac{\partial f_{3}}{\partial u_{3}}\;\frac{\partial f_{4}}{\partial u_{3}}\;\frac{\partial f_{5}}{\partial u_{3}}\;\frac{\partial f_{6}}{\partial u_{3}}
 \end{bmatrix}'$, at origin and $u=\begin{bmatrix}u_{1}\; u_{2}\; u_{3} \end{bmatrix}'$. The matrix $A$ is given by (\ref{A}) and $B=\begin{bmatrix}0\;1\;0\;0\;0\;0\\0\;0\;0\;1\;0\;0\\0\;0\;0\;0\;0\;1 \end{bmatrix}'$. The controllability matrix $Q$ is given by
 \begin{equation*}
 \begin{bmatrix}B\;AB\;A^{2}B\;A^{3}B\;A^{4}B\;A^{5}B \end{bmatrix}=\begin{bmatrix}
0 & 0 & 0 & 1 & 0 & 0 & 0 & 2 & 0 & -3 & 0 & 0 & 0 & -6 & 0 & 9 & 0 & 0\\
1 & 0 & 0 & 0 & 2 & 0 & -3 & 0 & 0 & 0 & -6 & 0 & 9 & 0 & 0 & 0 & 18 & 0\\
0 & 0 & 0 & 0 & 1 & 0 & -2 & 0 & 0 & 0 & -4 & 0 & 6 & 0 & 0 & 0 & 12 & 0\\
0 & 1 & 0 & -2 & 0 & 0 & 0 & -4 & 0 & 6 & 0 & 0 & 0 & 12 & 0 & -18 & 0 & 0\\
0 & 0 & 0 & 0 & 0 & 1 & 0 & 0 & 0 & 0 & 0 & 0 & 0 & 0 & 0 & 0 & 0 & 0\\
0 & 0 & 1 & 0 & 0 & 0 & 0 & 0 & 0 & 0 & 0 & 0 & 0 & 0 & 0 & 0 & 0 & 0
\end{bmatrix}.
\end{equation*} 
The rank of matrix $Q$ is 6, which is equal to the dimensions of the state $X\left(=6\right)$. Hence by the Kalman's condition the system is controllable if we add the thrusters in $r$, $\theta$ and $z$ direction. Figure-3 shows the trajectories of states of the system \eqref{eqm15} with initial state $[1\;2\;3\;4\;5\;6]'$ and desired final state  $[6\; 5\; 4\; 3\; 2\; 1]'$ respectively.
\begin{figure*}[ht]
\begin{tabular}{cc}
\includegraphics[width=9cm,height=9cm]{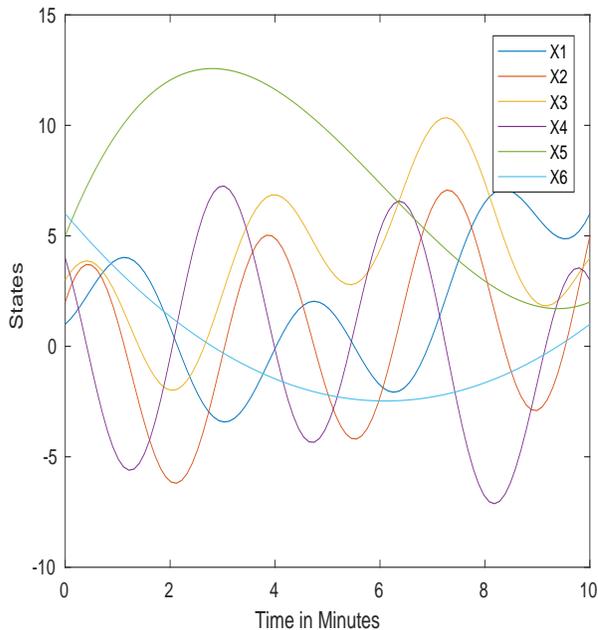}  \end{tabular}
\caption{State Control of the System under the effect of zonal harmonic $J_2$}
\end{figure*}
We can see from figure-3 that the initial state is steered to final state during the time interval $[0,10]$. The graph of the controllers i.e. thrusters in all the three directions $r$, $\theta$ and $z$ are shown in figure-4.
\begin{figure*}[ht]
\begin{tabular}{cc}
\includegraphics[width=9cm,height=9cm]{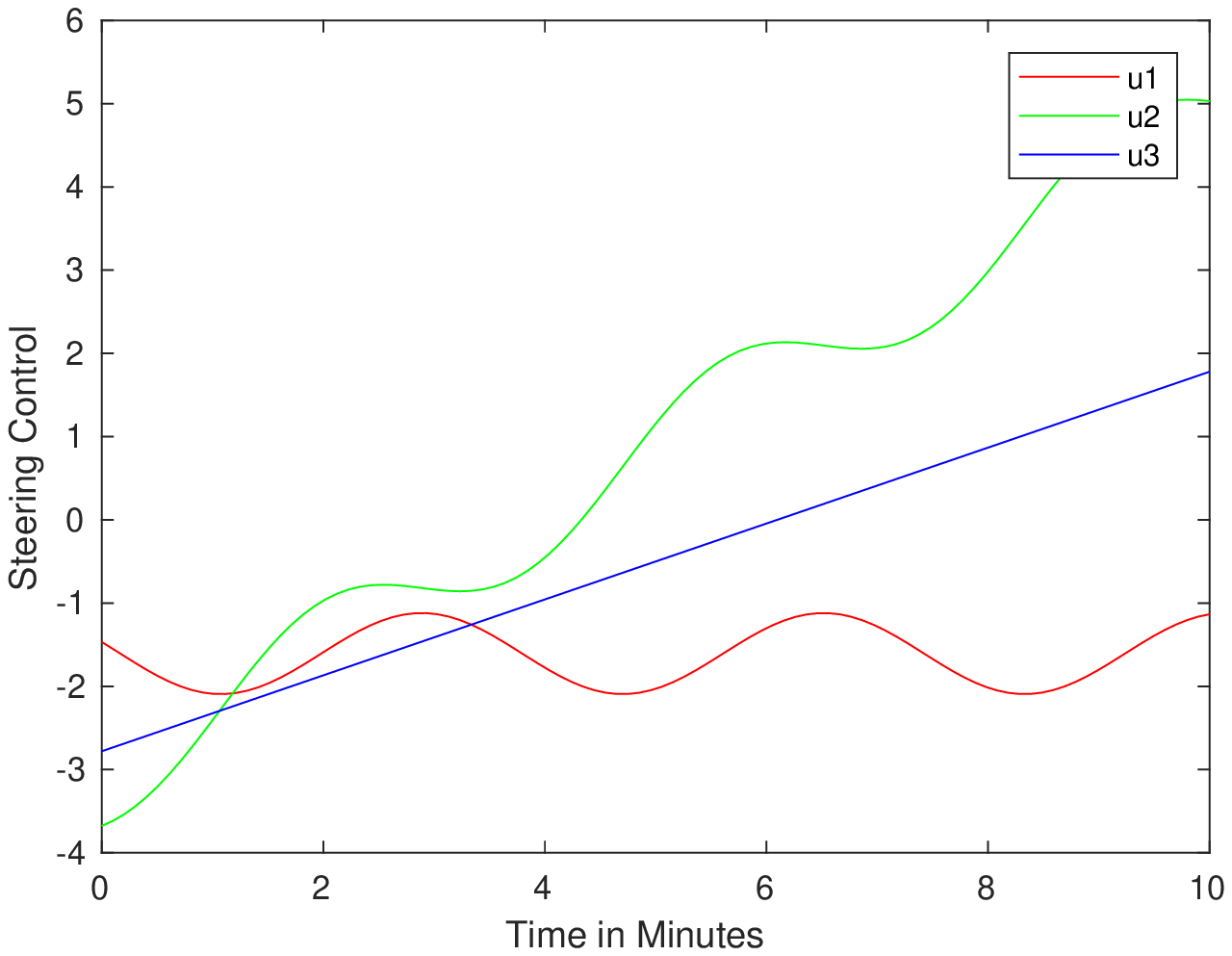}  
\end{tabular}
\caption{Steering Control of the System under the effect of zonal harmonic $J_2$}
\end{figure*}

\section{Conclusion}
\label{sec:4}
We have studied controllability analysis for seven different cases by applying controllers in (1) $r$- direction, (2) $\theta$- direction, (3) $z$- direction, (4) $r$ and $\theta$ directions, (5) $r$ and $z$ directions, (6) $\theta$ and $z$ directions and (7) $r$, $\theta$ and $z$ directions. Applying the Kalman's rank condition we found that, the system \eqref{eqm2} is uncontrollable if we apply thrusters i.e controllers in (1) $r$- direction, (2) $\theta$- direction, (3) $z$- direction, (4) $r$ and $\theta$ directions, (5) $r$ and $z$ directions, and it is controllable if thrusters are applied in (6) $\theta$ and $z$ directions and (7) $r$, $\theta$ and $z$ directions. 

From this study we found that to control the motion of the satellite under the effect of zonal harmonic $J_2$ we need to plug the controllers in the form of thrusters in all three directions. If the thruster in $r$ direction fails then also motion of satellite is controllable, but if thruster in any other direction(s) fail then the motion of satellite will become uncontrollable and it may hit the Earth's surface.


\begin{thebibliography}{}
\bibitem{Awange} J. L. Awange, E. W. Grafarend, B. Palancz, P. Zaletnyik, Algebraic Geodesy and Geoinformatics. Springer, Heidelberg, New York, (2010).
\bibitem{Palancz} B. Palancz, Application of Dixon resultant to satellite trajectory control by pole placement, Journal of Symbolic Computation, \textbf{50}, (2013), 79.
\bibitem{Raj} J. X. Raj, Analytical and Numerical Predictions for Near Earth's Satellite Orbits with KS Uniform Regular
Canonical Equations, PhD Thesis, Vikram Sarabhai Space Centre, India, (2007).
\bibitem{SS}  E. L. Stiefel and G. Scheifele, Linear and Regular Celestial Mechanics, Springer-Verlag, Berlin, Heidelberg, New York, (1971).
\bibitem{King} D. G. King-Hele, The effect of Earth's oblateness on the orbit of a near satellite, Proc. R. Soc. London A,
Math. Phys. Sci., \textbf{247}, (1958), 49.
\bibitem{Sehnal} L. Sehnal,The Earth upper atmosphere and the motion of Artificial Satellites,Publications of the Department of Astronomy- Beograd, \textbf{10}, (1980), 5.
\bibitem{KPTN} S. H. Knowles, J. E. Picon, S.E. Thonnard and A. C. Nicholas,The effect of Atmospheric drag on Satellite orbits during the Bastille day event, Solar Physics, \textbf{204}, (2001), 387.
\bibitem{YK} Q. Yan and V. Kapila, Analysis and Control of Satellite orbits around oblate Earth using perturbation method, Proc. $40^{th}$ IEEE Conf. Dec. Cont., Orlando, Florida, (2001), 1517.
\bibitem{Khalil} KH. I. Khalil,The drag exerted by an oblate rotating atmosphere on artificial satellite, Appl. Math. Mech., \textbf{23}, (2002), 1016.
\bibitem{BV} A. Bezdv\v{e}k and D. Vokrouhlick\'{y}, Semianalytic theory of motion for close-Earth spherical satellite including drag and gravitational perturbations, Planet. Spa. Sci., \textbf{52}, (2004), 1233.
\bibitem{HHB} I. A. Hassan, Z. M. Hayman and M. A. F. Basha, Pre-solution of perturbed motion of artificial satellite, Proc. First Middle East Africa IAU- Regional Meet.-1, (2008), 16.
\bibitem{CJ} W. Chen and W. Jing, Dyanamic equtions of relative motion around an oblate earth with air drag, Journal of Aero-space Engineering, \textbf{25}, (2012), 21.
\bibitem{RM} T. Reid and A. K. Misra, Formation flight of satellite in the presence of atmospheric drag, J. Aero. Engin. Sci. Appl., \textbf{3}, (2011), 64.
\bibitem{LSSC} D. Lee, J. C. Springmann, S.C. Spangelo and J. W. Cutler, Satellite dynamics simulator development using Lie group variational integrator, Proc. AIAA Mod. Sim. Tech., (2010), 1.
\bibitem{XTCY} G. Xu, X. Tianhe, W. Chen and T. Yeh, Analytical solution of satellite orbit disturbed by atmospheric drag, Mon. Not. R. Astron. Soc., \textbf{410}, (2011), 654.
\bibitem{AAAB} M. J. F. Al-Bermani, Abed Al-Ameer H. Ali, A. M. Al-Hashmi, A. S. Baron , Effect of atmospheric drag and zonal harmonic on Cosmos1484 satellite orbit, J. Kufa - Phys., \textbf{4}, (2012), 1.
\bibitem{Delhaise} F. Delhaise, Analytical treatment of air drag and earth oblateness effect upon an artificial satellite, Cel. Mecha. Dyna. Astron., \textbf{52}, (1991), 85.
\bibitem{Battin} R. H. Battin, An Introduction to Mathematics and Methods of Astrodynamics, AIAA Education Series, New York, (1987).
\bibitem{Wiesel} W. E. Wiesel, Modern Astrophysics, Aphelion Press, (2003).
\bibitem{Vallado} D. A. Vallado, Fundamentals of Astrodynamics and Applications, Microcosm Press and Kluwer Academic Publisher, (2004).
\bibitem{GWA} M. Grewal, L. Weill and A. Andrews, Global Positioning System, Intertial, Navigation and Integration, A John Wiley and Sons Inc. Publication, (2007).
\bibitem{SRPSCB} J. P. Sharma, B. S. Ratanpal, U. M. Pirzada, V. S. Shah, A. B. Chavda, N. B. Dave, Study of Effect of Perturbation Due ot Oblatness of Earth on Satellite, Proceeing of 19th Anual cu 4th International Conference of Gwalior Academy of Mathematical Sciences, Sardar Vallabhbhai National Institute of Technology, Surat, (2014) 338.
\bibitem{SRPS} J. P. Sharma, B. S. Ratanpal, U. M. Pirzada, V. S. Shah, Motion of Satellite under the Effect of Oblateness of
Earth and Atmospheric Drag, The Int. J. Analy. Exp. Model. Analy., \textbf{XI}, (2019), 2514. 
\bibitem{Hajo} B. B. Hajovsky, Satellite Formation Control by Atmospheric Drag, Master thesis, Airforce Institute of Technology, (2007).
\bibitem{Neo} K. Neokleous, Modelling and Control of Satellite's Geostationary Orbit, Diploma thesis, Lulea University of Technology, (2007).
\bibitem{Palancz1} B. Palancz, Numeric- Symbolic solution for satellite Trajectory Control for Pole Placement, Periodica Polytechnica, \textbf{57} (2013), 21.
\bibitem{Lamba} S. Lamba, Controllability, Observability and Stability of Artificial Satellite Problem, Master Thesis, National Institute of Technology, Jamshedpur, (2017).
\bibitem{Brockett} R. W. Brockett, Finite dimensional Linear Systems, John Wiley and Sons, Newyork, (1970).
\bibitem{Humi} M. Humi, $J_2$ Effect in Cylindrical Coordinates, J. Guidance, Contron and Dynamics, \textbf{30}, (2007), 263.
\end{thebibliography}
\end{document}